\documentclass{amsart}
\usepackage{hyperref}

\numberwithin{equation}{section}

\def\a{\alpha}
\def\b{\beta}
\def\d{\delta}
\def\D{\Delta}
\def\e{\varepsilon}
\def\f{\varphi}
\def\F{\Phi}
\def\g{\gamma}

\def\k{\kappa}

\def\L{\Lambda}

\def\z{\zeta}

\def\ov{\overline}

\def\SZ{Szeg\H o}
\def\({\left(}
\def\){\right)}
\def\[{\left[}
\def\]{\right]}

\font\tenopen = cmbx10 \font\sevenopen = cmbx7 \font\fiveopen =
cmbx5
\newfam\openfam

\textfont\openfam = \tenopen \scriptfont\openfam = \sevenopen
\scriptscriptfont\openfam = \fiveopen

\begin{document}

\title[Schur flows]{Schur flows and orthogonal polynomials on the unit circle}

\author{Leonid Golinskii}
\address{Mathematics Division,
Institute for Low Temperature Physics and Engineering\\
47\\ Lenin Ave\\ Kharkov, 61103 \\ Ukraine}
\email{\href{mailto:golinskii@ilt.kharkov.ua}{golinskii@ilt.kharkov.ua}}


\begin{abstract} The relation between the Toda lattices and
similar nonlinear chains and orthogonal polynomials on the real
line has been elaborated immensely for the last decades. We
examine another system of differential-difference equations known
as the Schur flow, within the framework of the theory of
orthogonal polynomials on the unit circle. This system can be
displayed in equivalent form as the Lax equation, and the
corresponding spectral measure undergoes a simple transformation.
The general result is illustrated on the modified Bessel measures
on the unit circle and the long time behavior of their Verblunsky
coefficients
\end{abstract}

\maketitle

\section{Introduction}

In 1975 J. Moser \cite{M1, M2} suggested a method for solution of
the finite Toda lattice equations (specifically, the Cauchy
problem for such lattices) based on the spectral theory of finite
Jacobi matrices. Later on Yu. M. Berezanskii \cite{Be1} adapted
this method to semi-infinite Toda lattices
\begin{eqnarray} \label{1.1}
 \nonumber 
  a_n' &=& a_n(b_{n+1}-b_n), \\
  b_n' &=& 2(a_n^2-a_{n-1}^2), \qquad
  n\in\mathbb{Z}_+=\{0,1,\ldots\}, \quad a_{-1}=0,
\end{eqnarray}
where $'$ means differentiation with respect to $t$, in the class
of bounded real $b$'s and positive $a$'s with the initial data
$\{b_n(0)=\overline{b_n(0)}, \ a_n(0)>0\}$. The key idea is to
compose a semi-infinite Jacobi matrix
\begin{equation}    \label{1.3}
    J=J(\{a_n\},\{b_n\})=
    \left(\begin{array}{ccccc}
                b_0(t)&a_0(t)&&&\\
                a_0(t)&b_1(t)&a_1(t)&&\\
                &a_1(t)&b_2(t)&a_2(t)&\\
                &&\ddots&\ddots&\ddots\\
            \end{array}\right)
\end{equation}
and trace the evolution of the matrix-valued function $J=J(t)$ and
its spectral characteristics. It turned out that (\ref{1.1}) can
be paraphrased in equivalent form in terms of $J$ itself (the Lax
equation)
\begin{eqnarray}
J'(t)&=&[A,J]=AJ-JA, \\
   A &=& \left(\begin{array}{ccccc}
                0&a_0(t)&&&\\
                -a_0(t)&0&a_1(t)&&\\
                &-a_1(t)&0&a_2(t)&\\
                &&\ddots&\ddots&\ddots\\
            \end{array}\right) = \pi(J):=J_+-J_-
\end{eqnarray}
with the standard notation $X_{\pm}$ for the upper (lower)
projection of a matrix $X$, as well as the corresponding spectral
(orthogonality) measure $d\mu(x,t)$ which undergoes a simple
transformation
\begin{equation}\label{1.5}
    d\mu(x,t)=e^{-xt}d\mu(x,0).
\end{equation}
Hence the solution of (\ref{1.1}) boils down to a combination of
the direct spectral problem (from $\{a_n(0), b_n(0)\}$ to
$d\mu(x,0)$) at $t=0$, plus (\ref{1.5}), plus the inverse spectral
problem (from $d\mu(x,t)$ to $\{a_n(t), b_n(t)\}$) at $t>0$.

The theory of orthogonal polynomials on the real line plays one of
the first fiddles in the performance (albeit not entering the
final result directly). For instance, it furnishes a nice setting
for solving the inverse spectral problem. There is a parallel
theory of orthogonal polynomials on the unit circle (OPUC) which
has experienced a splash of activity lately thanks to primarily
Simon's disquisition \cite{S1, S2}. So the question arises
naturally whether there exist nonlinear chains (so to say, the
``Toda lattices for the unit circle'') which can be handled by the
similar method. The main goal of the present paper is to develop
the ``Moser--Berezanskii scheme for the unit circle'' based on the
spectral theory of a certain class of unitary matrices in
application to a system of nonlinear differential-difference
equations known as the {\it Schur flow}.

We begin with some basics on orthogonal polynomials on the unit
circle (OPUC). Given a nontrivial (i.e., not a finite combination
of delta functions) probability measure $\mu$ on the unit circle
$\mathbb{T}$ with the moments
 $$ \mu_k:=\int_{\mathbb{T}} \z^{-k}d\mu, \qquad
 k\in\mathbb{Z}=\{0,\pm1,\ldots\}, $$
we define the monic orthogonal polynomials $\F_n(z,\mu)$ (or just
$\F_n$ if $\mu$ is understood) by
 \begin{equation} \label{1.6}
\int_{\mathbb{T}} \z^{-k}\F_n(\z)d\mu=0, \ \ k=0,\ldots, n-1;
\quad \F_n(z)=z^n+l_nz^{n-1}+\ldots+\F_n(0).
\end{equation}
Clearly such system is uniquely determined and
 \begin{equation}\label{1.7}
\int_{\mathbb{T}} \F_m(\z)\ov{\F_n(\z)}\,d\mu=0, \qquad m\not= n.
\end{equation}
The orthonormal polynomials $\f_n=\k_n\F_n$, $\k_n>0$ enjoy the
property
$$ \int_{\mathbb{T}} \f_m(\z)\ov{\f_n(\z)}\,d\mu=\delta_{mn}. $$

 A key role throughout the whole OPUC
theory is played by the sequences of complex numbers
$\{\a_n\}_{n\geq0}$, $|\a_n|<1$,
\begin{equation}\label{1.8}
   \a_n=\a_n(\mu)=-\ov{\F_{n+1}(0)}, \quad
 n\in\mathbb{Z}_+, \quad \a_{-1}:=-1,
\end{equation}
known as the {\it Verblunsky coefficients} or parameters of OPUC
system. Firstly, due to the celebrated Verblunsky theorem, there
is one-one correspondence between the class $\mathcal{P}$ of all
nontrivial probability measures on $\mathbb{T}$ and the set
$\mathbb{D}^\infty$, so each sequence of complex numbers
$\{\gamma_n\}_{n\geq0}$ from the open unit disk $\mathbb{D}$ comes
up as a system of parameters for uniquely determined measure
$\mu\in\mathcal{P}$. Secondly, Verblunsky coefficients (\ref{1.8})
enter the Szeg\H o recurrence relations given in the vector form
by
\begin{equation}\label{1.10}
\left[%
\begin{array}{c}
  \F_{n+1}(z) \\
  \F_{n+1}^*(z) \\
\end{array}%
\right] = T_n(z)
\left[%
\begin{array}{c}
  \F_n(z) \\
  \F_n^*(z) \\
\end{array}%
\right], \qquad T_n(z)=
\left(%
\begin{array}{cc}
  z & -\ov\a_n \\
  -\a_n z & 1 \\
\end{array}%
\right)
\end{equation}
is called the \SZ\ matrix, and so
$$
\left[%
\begin{array}{c}
  \F_{n+1}(z) \\
  \F_{n+1}^*(z) \\
\end{array}%
\right] = \mathcal{T}_n(z)
\left[%
\begin{array}{c}
  1 \\
  1 \\
\end{array}%
\right], \qquad \mathcal{T}_n(z)= T_n(z)T_{n-1}(z)\ldots T_0(z)
$$
is the transfer matrix. So both monic orthogonal and orthonormal
polynomials are completely determined by the sequence $\a_n$, the
latter because of the equality
$$ \k_n^{-2}=\prod_{j=0}^{n-1} (1-|\a_j|^2). $$

To complete with the basic properties of OPUC let us mention
explicit determinant formulae for both monic polynomials and
Verblunsky coefficients in terms of the moments of the
orthogonality measure:
$$\Phi_n(z)=\frac1{D_n} \left|
\begin{array}{cccc}
  \mu_0 & \mu_{-1} & \ldots & \mu_{-n} \\
  \mu_1 & \mu_0 & \ldots & \mu_{-n+1} \\
  \vdots & \vdots & \vdots & \vdots \\
  \mu_{n-1} & \mu_{n-2} & \ldots & \mu_{-1} \\
  1 & z & \ldots & z^n \\
\end{array}%
\right|, \quad D_{n+1}:=\det\|\mu_{k-j}\|_{k,j=0}^n,$$
\begin{equation}\label{1.12}
 \F_n(0)=-\ov\a_{n-1}=\frac{(-1)^n}{D_n} \left|
\begin{array}{ccc}
  \mu_{-1} & \ldots & \mu_{-n} \\
  \mu_0 & \ldots & \mu_{-n+1} \\
  \vdots & \vdots & \vdots \\
  \mu_{n-2} & \ldots & \mu_{-1} \\
\end{array}%
\right|.
 \end{equation}

One of the most interesting developments in the theory of OPUC in
recent years is the discovery by Cantero, Moral, and Vel\'azquez
\cite{CMV} of a matrix realization for multiplication by $\z$ on
$L^2(\mathbb{T}, \mu)$ which is of finite band size (i.e.,
$|\langle\z\chi_m,\chi_n\rangle|=0$ if $|m-n|>k$ for some $k$; in
this case, $k=2$ to be compared with $k=1$ for the real line
case). Their basis (complete, orthonormal system) $\{\chi_n\}$ is
obtained by orthonormalizing the sequence
$1,\z,\z^{-1},\z^2,\z^{-2},\ldots$. Remarkably, the $\chi$'s can
be expressed in terms of $\f$'s and $\f^*$'s
 (see \cite [Proposition 4.2.2]{S1})
\begin{equation}\label{1.14}
 \chi_{2n}(z)=z^{-n}\f_{2n}^*(z),\qquad
 \chi_{2n+1}(z)=z^{-n}\f_{2n+1}(z), \quad n\in\mathbb{Z}_+,
\end{equation}
and the matrix elements
 $$ \mathcal{C}(\mu)=\|c_{nm}\|=\langle\z\chi_m,\chi_n\rangle,
    \qquad m,n\in\mathbb{Z}_+  $$
in terms of Verblunsky coefficients
\begin{equation}\label{1.17}
\mathcal{C}(\{\a_n\})=\left(%
\begin{array}{cccccc}
  \ov\a_0 & \ov\a_1\rho_0 & \rho_0\rho_1 & 0 & 0 & \ldots \\
  \rho_0 & -\ov\a_1\a_0 & -\a_0\rho_1 & 0 & 0 & \ldots \\
  0 & \ov\a_2\rho_1 & -\ov\a_2\a_1 & \ov\a_3\rho_2 & \rho_2\rho_3 & \ldots \\
  0 &\rho_1\rho_2 & -\a_1\rho_2 & -\ov\a_3\a_2 & -\a_2\rho_3 & \ldots \\
  0 & 0 & 0 & \ov\a_4\rho_3 & -\ov\a_4\a_3 & \ldots \\
  \vdots & \vdots & \vdots & \vdots & \vdots & \vdots \\
\end{array}%
\right) \end{equation}
  with $\rho_n^2:=1-|\a_n|^2$, $0<\rho\leq 1$.

There is an important relation between CMV matrices and monic
orthogonal polynomials akin to the well-known property of
orthogonal polynomials on the real line:
\begin{equation}\label{1.18}
 \F_n(z)=\det(zI_n-\mathcal{C}^{(n)}),
\end{equation}
where $\mathcal{C}^{(n)}$ is the principal $n\times n$ block of
$\mathcal{C}$ (see, e.g., \cite [p. 271, formula (4.2.47)]{S1}).

 The CMV matrices $\mathcal{C}$ play much the same role in our
study of the Schur flows that Jacobi matrices (\ref{1.3}) in the
Toda lattices setting.

We are in a position now to announce our main result.

\noindent
{\bf Theorem 1}.  {\it Let $\alpha_n(t)$,
$n\in\mathbb{Z}_+$, be a sequence of complex valued functions with
$|\alpha_n(t)|<1$ for $t\geq 0$ and let $\alpha_{-1}=-1$. The
following three statements are equivalent.
\begin{enumerate}
    \item $\alpha_n$ solve the Schur flow equations
    \begin{equation}\label{1.20}
    \a_n'(t)=(1-|\a_n(t)|^2)(\a_{n+1}(t)-\a_{n-1}(t)),
    \qquad t>0;
\end{equation}
    \item The CMV matrices $\mathcal{C}(t)$ satisfy the Lax
    equation
    \begin{equation}\label{1.21}
    \mathcal{C}'(t)=[A,\mathcal{C}],
\end{equation}
where $A$ is an upper triangular and tridiagonal matrix
\begin{equation}\label{1.22}
A=\|a_{mn}\|_0^\infty=\left(%
\begin{array}{cccccc}
  \Re\ov\a_0 & \rho_0\ov\D_0 & \rho_0\rho_1 & {} & {} & {} \\
  {} & -\Re\ov\a_1\a_0 & \rho_1\D_1 & \rho_1\rho_2 & {} & {} \\
  {} & {} & -\Re\ov\a_2\a_1 & \rho_2\ov\D_2 & \rho_2\rho_3 & {} \\
  {} & {} & {} & \ddots & \ddots & \ddots \\
\end{array}%
\right),  \end{equation} $\D_n:=\a_{n+1}(t)-\a_{n-1}(t)$;
    \item The orthogonality measure $\mu(\z,t)$, having $\a_n(t)$
    as its Verblunsky coefficients, satisfies
    \begin{equation}\label{1.23}
    d\mu(\zeta,t)=C(t)e^{t(\zeta+\zeta^{-1})}\,d\mu(\zeta,0),
\end{equation}
where $C$ is a normalizing factor.
\end{enumerate} }

 We refer to (\ref{1.23}) as the {\it Bessel transformation} of
 the measure $d\mu=d\mu(\z,0)$.

\smallskip\noindent
{\sl Remarks}. 1. It is not hard to see that if $\{\a_n(t)\}$
solves (\ref{1.20}) with the initial data $|\a_n(0)|<1$, then
$|\a_n(t)|<1$ for each $t>0$. Indeed, assume for the contrary that
for some $n\geq 0$ there is $t_0>0$ such that $|\a_n(t_0)|=1$ and
$|\a_n(t)|<1$ for all $0<t<t_0$. It follows from (\ref{1.20}) and
$\rho_n^2=1-|\a_n|^2$ that
$$ (\rho_n^2)'=2\rho_n\rho_n'=-\a_n'\ov\a_n-\a_n\ov\a_n'
   =-2\rho_n^2\Re(\ov\a_{n+1}\a_n-\ov\a_n\a_{n-1}), $$
and so
 \begin{equation}\label{1.30}
 \rho_n'=-\rho_n\Re(\ov\a_{n+1}\a_n-\ov\a_n\a_{n-1}), \quad
 0<t<t_0.
 \end{equation}
Hence
 \begin{equation}\label{1.31}
 \rho_n(t)=\rho_n(0)\exp\left\{-\int_0^t
\Re(\ov\a_{n+1}(s)\a_n(s)-\ov\a_n(s)\a_{n-1}(s))ds\right\},
\end{equation}
and the right hand side is bounded away from zero as $t\to t_0$,
whereas the left hand side tends to zero. The contradiction shows
that $|\a_n(t)|<1$ for all $n\in\mathbb{Z}_+$ and $t>0$.
Therefore, by (3) the solution of the initial-boundary-value
(IBV) problem
    \begin{equation}\label{1.33}
    \a_n'(t)=(1-|\a_n(t)|^2)(\a_{n+1}(t)-\a_{n-1}(t)),
    \quad t>0; \qquad |\a_n(0)|<1,
\end{equation}
$n\in\mathbb{Z}_+$, $\a_{-1}=-1$ exists and unique.

2. We can modify the second statement by observing that
$$ A+A^*=\mathcal{C}+\mathcal{C}^*=
\left(%
\begin{array}{cccccc}
  2\Re\ov\a_0 & \rho_0\ov\D_0 & \rho_0\rho_1 & {} & {} & {} \\
  \rho_0\D_0 & -2\Re\ov\a_1\a_0 & \rho_1\D_1 & \rho_1\rho_2 & {} & {} \\
  \rho_0\rho_1 & \rho_1\ov\D_1 & -2\Re\ov\a_2\a_1 & \rho_2\ov\D_2 &
  \rho_2\rho_3 & {} \\
  \vdots & \vdots & \vdots & \ddots & \ddots & \ddots \\
\end{array}%
\right), $$
and so
\begin{equation}\label{1.35}
\mathcal{C}'(t)=[B,\mathcal{C}],
\end{equation}
\begin{eqnarray}\label{1.37}
 \nonumber B = \frac{(\mathcal{C}+\mathcal{C}^*)_+
-(\mathcal{C}+\mathcal{C}^*)_-}2 &=&
\frac12 \left(%
\begin{array}{cccccc}
  0 & \rho_0\ov\D_0 & \rho_0\rho_1 & {} & {} & {} \\
  -\rho_0\D_0 & 0 & \rho_1\D_1 & \rho_1\rho_2 & {} & {} \\
  -\rho_0\rho_1 & -\rho_1\ov\D_1 & 0 & \rho_2\ov\D_2 & \rho_2\rho_3 & {} \\
  \vdots & \vdots & \vdots & \ddots & \ddots & \ddots \\
\end{array}%
\right) \\
 {} &=& A-\frac{\mathcal{C}+\mathcal{C}^*}2=-B^*,
\end{eqnarray}
 which makes it closer to its counerpart in the Toda lattices setting.

So, once again, the solution of IBV problem (\ref{1.33}) amounts
to a combination of the direct and inverse spectral problems with
(\ref{1.23}) in between. Note that the orthogonality measure
$\mu(\z,0)$ can be retrieved from  the initial data $\a_n(0)$ by
either the Spectral Theorem for the CMV matrix $\mathcal{C}(0)$
(\ref{1.17}) or via orthonormal polynomials, since $\mu(\z,0)$
arises as a *-weak limit of the sequence of measures
$|\f_n|^{-2}dm$, $dm$ being a normalized Lebesgue measure on
$\mathbb{T}$ (Rakhmanov's theorem). In turn, the Verblunsky
coefficients $\a_n(t)$ are recovered from the measure $\mu(\z,t)$
by (\ref{1.12}).

The Schur flow (\ref{1.20}) emerged in \cite{AL1, AL2} under the
name {\it discrete modified KdV equation}, as a spatial
discretization of the modified Korteweg--de Vries equation
$$ \partial_t f=6f^2\partial_x f-\partial_x^3 f. $$
In \cite{FG} the authors deal with finite real Schur flows and
suggest two more distinct Lax equations based on the Hessenberg
matrix representation of the multiplication operator (see also
\cite{AG}). In \cite{MN1, MN2} the Bessel modification of measures
appeared and a part of our main result which concerns
$(3)\Rightarrow (1)$ is proved. In a recent paper \cite{Ne1} the
author deals with the Poisson structure and Lax pairs for the
Ablowits--Ladik systems closely related to the Schur flows. The
latter can also be viewed as the zero-curvature equation for the
Szeg\H o matrices (cf. \cite{GGH})
$$  T_n'(z,t) + T_n(z,t) W_n(z,t) -W_{n+1}(z,t) T_n(z,t)=0, $$
$$ W_n(z,t):=\left(%
\begin{array}{cc}
  z+1-\a_{n-1}\ov\a_n & -\ov\a_n-\ov\a_{n-1}z^{-1} \\
  -\a_{n-1}z-\a_n & 1-\ov\a_{n-1}\a_n+z^{-1} \\
\end{array}%
\right). $$

We proceed as follows. The proof of our main result is presented
In Sections 2 and 3 with some comments on the general IBV problem
and doubly infinite systems. We also show that some properties of
Verblunsky coefficients for the Bessel transformed measures, such
as the rate of decay, are inherited from those of the initial
data. In Section 4 we study the modified Bessel measures and
polynomials on the unit circle, a nice example which corresponds
to the zero initial conditions in our setting. In this case the
long time behavior of the Verblunsky coefficients can be obtained.

\section{Proof of Theorem 1: $(1)\Leftrightarrow (2)$.}

\noindent
$(1)\Rightarrow (2)$. Once the Lax pair is enunciated,
the proof goes through by brute force computation which is much
more involved compared to the Toda case.

Let $2\e_m:=1-(-1)^m$,\ $m\in\mathbb{Z}_+$, and $\e_{-1}=1$, so
$\{\e_m\}_{m\geq 0}=\{0,1,0,1,\ldots\}$,
$$ \e_m+\e_{m+1}=1, \qquad \e_m\e_{m+1}=0, \qquad \e_m-\e_{m+1}=(-1)^{m+1}. $$
It is instructive to write the diagonals of $\mathcal{C}$
(\ref{1.17}) in a unique way
\begin{eqnarray}\label{}
 c_{mm} &=& -\ov\a_m\a_{m-1},\\
 c_{m+2,m} &=& \rho_m\rho_{m+1}\e_m, \qquad
 c_{m,m+2}=\rho_m\rho_{m+1}\e_{m+1},
\end{eqnarray}
and
\begin{eqnarray}
  c_{m+1,m} &=& \ov\a_{m+1}\rho_m\e_m-\a_{m-1}\rho_m\e_{m+1}, \\
  c_{m,m+1} &=& \ov\a_{m+1}\rho_m\e_{m+1}-\a_{m-1}\rho_m\e_m.
\end{eqnarray}
In the same vein for the matrix entries of $A$ (\ref{1.22})
\begin{equation}\label{}
 a_{mm}=-\Re\ov\a_m\a_{m-1}, \qquad a_{m,m+2}=\rho_m\rho_{m+1},
\end{equation}
and
\begin{equation}\label{}
 a_{m,m+1}=\rho_m\ov\D_m\e_{m+1}+\rho_m\D_m\e_m.
\end{equation}

Next, it follows from (\ref{1.20}) and (\ref{1.30}) that
\begin{eqnarray}\label{}
\nonumber (\rho_m\rho_{m+1})'
 &=& -\rho_m\rho_{m+1}\Re(\ov\a_{m+2}\a_{m+1}-\ov\a_m\a_{m-1}),\\
\nonumber (\ov\a_m\a_{m-1})' &=&
 \a_{m-1}\rho_m^2\ov\D_m+\ov\a_m\rho_{m-1}^2\D_{m-1},
\end{eqnarray}
and
\begin{eqnarray}
\nonumber  (\a_{m-1}\rho_{m})' &=&
\[\rho_{m-1}^2\D_{m-1}-\a_{m-1}
  \Re(\ov\a_{m+1}\a_{m}-\ov\a_{m}\a_{m-1})\]\rho_m, \\
\nonumber  (\ov\a_{m+1}\rho_m)' &=&
  \[\rho_{m+1}^2\ov\D_{m+1}-\ov\a_{m+1}\Re(\ov\a_{m+1}\a_m-\ov\a_m\a_{m-1})\]\rho_m.
\end{eqnarray}
 Hence, for derivatives of the CMV matrix entries we have now
\begin{eqnarray}\label{}
 c_{mm}' &=& -(\ov\a_m\a_{m-1})'
 =-\a_{m-1}\rho_m^2\ov\D_m-\ov\a_{m}\rho_{m-1}^2\D_{m-1},\\
 c_{m,m+2}' &=& -\rho_m\rho_{m+1}
 \Re(\ov\a_{m+2}\a_{m+1}-\ov\a_m\a_{m-1})\e_{m+1}, \\
 c_{m,m-2}' &=& -\rho_{m-2}\rho_{m-1}
 \Re(\ov\a_m\a_{m-1}-\ov\a_{m-2}\a_{m-3})\e_m,
\end{eqnarray}
and
\begin{eqnarray}
 \nonumber \ \ \ \ \ \ \ \ c_{m,m+1}' &=& \[\rho_{m+1}^2\ov\D_{m+1}-\ov\a_{m+1}
  \Re(\ov\a_{m+1}\a_m-\ov\a_m\a_{m-1})\]\rho_m\e_{m+1} \\
  {} &-& \[\rho_{m-1}^2\D_{m-1}-\a_{m-1}
  \Re(\ov\a_{m+1}\a_m-\ov\a_m\a_{m-1})\]\rho_m\e_{m},
\end{eqnarray}
\begin{eqnarray}
\nonumber  \ \ \ \ \ \ \ \ c_{m,m-1}' &=&
\[\rho_{m}^2\ov\D_m-\ov\a_{m}
  \Re(\ov\a_{m}\a_{m-1}-\ov\a_{m-1}\a_{m-2})\]\rho_{m-1}\e_{m+1} \\
{} &-& \[\rho_{m-2}^2\D_{m-2}-\a_{m-2}
  \Re(\ov\a_{m}\a_{m-1}-\ov\a_{m-1}\a_{m-2})\]\rho_{m-1}\e_{m}.
\end{eqnarray}

Let $K=\|K_{mn}\|=[A,\mathcal{C}]$,
$K_{mn}=\sum_i(a_{mi}c_{in}-c_{mi}a_{in})$. Since both $A$ and
$\mathcal{C}$ are of band size 2, i.e., $a_{mn}=c_{mn}=0$ for
$|m-n|>2$, and $A$ is an upper triangular, i.e., $a_{mn}=0$ for
$m>n$, we actually have
\begin{eqnarray}\label{2.16}
 K_{mn} &=& a_{mm}c_{mn}+a_{m,m+1}c_{m+1,n}+a_{m,m+2}c_{m+2,n} \\
 \nonumber &-& c_{mn}a_{nn}-c_{m,n-1}a_{n-1,n}-c_{m,n-2}a_{n-2,n}.
\end{eqnarray}

We want to show that
\begin{equation}\label{2.17}
 c_{m,m+j}'(t)=K_{m,m+j}(t), \qquad m,m+j\in\mathbb{Z}_+.
\end{equation}
To this end we will plug (2.1)--(2.6) into (\ref{2.16}) and
compare the outcome with (2.7)--(2.11).
 \footnote{One needs some paper and patience to slug one's way
 through the lengthy calculation.}
For $j\geq 5$ and $j\leq -3$ the equality holds for trivial
reason, as the both sides in (\ref{2.17}) vanish. For $j=4$
\begin{eqnarray*}
  K_{m,m+4} &=& a_{m,m+2}c_{m+2,m+4}-c_{m,m+2}a_{m+2,m+4} \\
  {}  &=& \rho_{m}\rho_{m+1}\rho_{m+2}\rho_{m+3}(\e_{m+3}-\e_{m+1})=0.
\end{eqnarray*}
For $j=3$ by (2.2) and (2.4)--(2.6)
\begin{eqnarray*}
  K_{m,m+3} &=& a_{m,m+1}c_{m+1,m+3}+a_{m,m+2}c_{m+2,m+3} \\
  &-& c_{m,m+1}a_{m+1,m+3}-c_{m,m+2}a_{m+2,m+3}= \end{eqnarray*}
$$ \rho_m\rho_{m+1}\rho_{m+2}\[\e_m(\D_m-\a_{m+1}+\a_{m-1})+
\e_{m+1}(\ov\a_{m+3}-\ov\a_{m+1}-\ov\D_{m+2})\]=0,$$
that is consistent with the banded structure of size 2 of
$\mathcal{C}'$.

The main work begins when $|j|=0,1,2$.

\noindent
 1. $j=-2$. We have by (2.2) and (2.5)
\begin{eqnarray}\label{2.20}
  K_{m,m-2} &=& c_{m,m-2}(a_{mm}-a_{m-2,m-2}) \\
  {} &=&
  \nonumber \rho_{m-2}\rho_{m-1}\Re(\ov\a_{m-2}\a_{m-3}-\ov\a_m\a_{m-1})\e_m,
\end{eqnarray}
and so (\ref{2.17}) holds by (2.9).

 \noindent
 2. $j=-1$. Write $K_{m,m-1}=K^{(1)}_{m,m-1}+K^{(2)}_{m,m-1}$ with
\begin{eqnarray*}
  K^{(1)}_{m,m-1} &=& c_{m,m-1}(a_{mm}-a_{m-1,m-1}) \\
  {} &=&
  (\ov\a_m\rho_{m-1}\e_{m+1}-\a_{m-2}\rho_{m-1}\e_m)
  \Re(\ov\a_{m-1}\a_{m-2}-\ov\a_m\a_{m-1}),\\
  K^{(2)}_{m,m-1} &=& a_{m,m+1}c_{m+1,m-1}-a_{m-2,m-1}c_{m,m-2} \\
  {} &=&
  \rho_{m-1}(\rho_m^2\ov\D_m\e_{m+1}-\rho_{m-2}^2\D_{m-2}\e_m)
\end{eqnarray*}
and hence
\begin{eqnarray}\label{2.21}
  \ \ \ \ \ \ \ \ K_{m,m-1} &=&
\[\rho_{m}^2\ov\D_m-\ov\a_{m}
  \Re(\ov\a_{m}\a_{m-1}-\ov\a_{m-1}\a_{m-2})\]\rho_{m-1}\e_{m+1} \\
\nonumber {} &-& \[\rho_{m-2}^2\D_{m-2}-\a_{m-2}
  \Re(\ov\a_{m}\a_{m-1}-\ov\a_{m-1}\a_{m-2})\]\rho_{m-1}\e_{m}.
\end{eqnarray}
Now (\ref{2.17}) follows from (2.11).

 \noindent
 3. $j=0$. Write $K_{mm}=K^{(1)}_{mm}+K^{(2)}_{mm}$ with
 \begin{eqnarray*}
K^{(1)}_{mm} &=&
a_{m,m+2}c_{m+2,m}-c_{m,m-2}a_{m-2,m}=\[(\rho_m\rho_{m+1})^2-
(\rho_{m-2}\rho_{m-1})^2\]\e_m, \\
  K^{(2)}_{mm} &=& a_{m,m+1}c_{m+1,m}-c_{m,m-1}a_{m-1,m} \\
  {} &=&
  (\ov\a_{m+1}\rho_m^2\D_m+\a_{m-2}\rho_{m-1}^2\ov\D_{m-1})\e_m \\
  &-& (\a_{m-1}\rho_m^2\ov\D_m+\ov\a_{m}\rho_{m-1}^2\D_{m-1})\e_{m+1}.
\end{eqnarray*}
But
$$ (\rho_m\rho_{m+1})^2 - (\rho_{m-2}\rho_{m-1})^2 +
          \ov\a_{m+1}\rho_m^2\D_m+\a_{m-2}\rho_{m-1}^2\ov\D_{m-1}
  =  $$
$$  \rho_m^2(1-\ov\a_{m+1}\a_{m-1})
  -\rho_{m-1}^2(1-\ov\a_{m}\a_{m-2})
     =  -\a_{m-1}\rho_m^2\ov\D_m-\ov\a_{m}\rho_{m-1}^2\D_{m-1}, $$
and so by (2.7)
$$
K_{mm}=-(\a_{m-1}\rho_m^2\ov\D_m+\ov\a_{m}\rho_{m-1}^2\D_{m-1})(\e_m+\e_{m+1})=c_{mm}'.
$$

 \noindent
 4. $j=1$. Now
$K_{m,m+1}=K^{(1)}_{m,m+1}+K^{(2)}_{m,m+1}+K^{(3)}_{m,m+1}$ with
\begin{eqnarray*}
  K^{(1)}_{m,m+1} &=& c_{m,m+1}(a_{mm}-a_{m+1,m+1}) \\
  {} &=&
  (\ov\a_{m+1}\e_{m+1}-\a_{m-1}\e_m)\rho_m\Re(\ov\a_{m+1}\a_m-\ov\a_m\a_{m-1}),\\
  K^{(2)}_{m,m+1} &=& a_{m,m+1}(c_{m+1,m+1}-c_{mm}) \\
  {} &=&
  -(\ov\D_m\e_{m+1}-\D_m\e_m)\rho_m(\ov\a_{m+1}\a_m-\ov\a_m\a_{m-1}),\\
  K^{(3)}_{m,m+1} &=& a_{m,m+2}c_{m+2,m+1}-a_{m-1,m+1}c_{m,m-1} \\
  {} &=&
  \[(\ov\a_{m+2}\rho_{m+1}^2-\ov\a_m\rho_{m-1}^2)\e_{m+1}-
  (\a_m\rho_{m+1}^2-\a_{m-2}\rho_{m-1}^2)\e_m\]\rho_m.
\end{eqnarray*}
Hence $\rho_m^{-1}K_{m,m+1}=u_m\e_m+v_m\e_{m+1}$ with
\begin{eqnarray*}
  u_m &=& -\a_{m-1}\Re(\ov\a_{m+1}\a_m-\ov\a_m\a_{m-1})+
  \a_{m-2}\rho_{m-1}^2-\a_m\rho_{m+1}^2 \\ &-&
  \D_m(\ov\a_{m+1}\a_m-\ov\a_m\a_{m-1}), \\
  v_m &=&
  \ov\a_{m+1}\Re(\ov\a_{m+1}\a_m-\ov\a_m\a_{m-1})+
  \ov\a_{m+2}\rho_{m+1}^2-\ov\a_m\rho_{m-1}^2 \\ &-&
  \ov\D_m(\ov\a_{m+1}\a_m-\ov\a_m\a_{m-1}).
\end{eqnarray*}
Next,
\begin{eqnarray*}
  u_m^{(1)} &:=& \a_{m-2}\rho_{m-1}^2-\a_m\rho_{m+1}^2+
  (\a_{m-1}-\a_{m+1})(\ov\a_{m+1}\a_m-\ov\a_m\a_{m-1}) \\
  {} &=& \a_{m+1}\ov\a_m\a_{m-1}-\ov\a_m\a_{m-1}^2+
  \ov\a_{m+1}\a_m\a_{m-1}-\a_m+\a_{m-2}\rho_{m-1}^2 \\
  {} &=& 2\a_{m-1}\Re\ov\a_{m+1}\a_m-2\a_{m-1}\Re\ov\a_m\a_{m-1}-
  (\a_m-\a_{m-2})\rho_{m-1}^2,
\end{eqnarray*}
and so
$$
u_m=\a_{m-1}\Re(\ov\a_{m+1}\a_m-\ov\a_m\a_{m-1})-\rho_{m-1}^2\D_{m-1}.
$$
In exactly the same way
$$
v_m=-\ov\a_{m+1}\Re(\ov\a_{m+1}\a_m-\ov\a_m\a_{m-1})+\rho_{m+1}^2\ov\D_{m+1},
$$
and finally
\begin{eqnarray}\label{2.22}
  \ \ \ \ \ \ \ \ K_{m,m+1} &=&
\[\rho_{m+1}^2\ov\D_{m+1}-\ov\a_{m+1}
  \Re(\ov\a_{m+1}\a_m-\ov\a_m\a_{m-1})\]\rho_m\e_{m+1} \\
\nonumber {} &-& \[\rho_{m-1}^2\D_{m-1}-\a_{m-1}
  \Re(\ov\a_{m+1}\a_m-\ov\a_m\a_{m-1})\]\rho_m\e_{m}.
\end{eqnarray}
We come to (\ref{2.17}) on account of (2.10).

\noindent
 5. $j=2$. Write $K_{m,m+2}=K^{(1)}_{m,m+2}+K^{(2)}_{m,m+2}+K^{(3)}_{m,m+2}$ with
\begin{eqnarray*}
  K^{(1)}_{m,m+2}\!\! &=&\!\! a_{m,m+2}(c_{m+2,m+2}-c_{mm})=
  \rho_m\rho_{m+1}(\ov\a_m\a_{m-1}-\ov\a_{m+2}\a_{m+1})(\e_m+\e_{m+1}), \\
  K^{(2)}_{m,m+2}\!\! &=&\!\! c_{m,m+2}(a_{mm}-a_{m+2,m+2})=
  \rho_m\rho_{m+1}\e_{m+1}\Re(-\ov\a_m\a_{m-1}+\ov\a_{m+2}\a_{m+1}),\\
  K^{(3)}_{m,m+2}\!\! &=&\!\! a_{m,m+1}c_{m+1,m+2}-c_{m,m+1}a_{m+1,m+2} \\
  &=&
  \rho_m\rho_{m+1}\[(\ov\a_{m+2}\D_m+\a_{m-1}\ov\D_{m+1})\e_m-
  (\a_m\ov\D_m+\ov\a_{m+1}\D_{m+1})\e_{m+1}\].
\end{eqnarray*}
We have
\begin{eqnarray*}
  \ov\a_{m+2}\D_m+\a_{m-1}\ov\D_{m+1} &=& -\ov\a_m\a_{m-1}+\ov\a_{m+2}\a_{m+1} \\
  \a_m\ov\D_m+\ov\a_{m+1}\D_{m+1} &=& -\a_m\ov\a_{m-1}+\a_{m+2}\ov\a_{m+1}
\end{eqnarray*}
and it follows from (2.8) that
\begin{equation}\label{2.23}
K_{m,m+2}=-\rho_m\rho_{m+1}\Re(\ov\a_{m+2}\a_{m+1}-\ov\a_m\a_{m-1})\e_{m+1},
\end{equation}
so (\ref{2.17}) holds again. The proof is complete.

\medskip
\noindent
$(2)\Rightarrow (1)$. The problem we are faced with is
that, in contrast to the Toda lattices, no $\a$'s in a pure form
appear among the matrix entries of $\mathcal{C}$.

Write
$|\ov\a_{n+1}\rho_n|^2=(1-\rho_{n+1}^2)\rho_n^2=\rho_n^2-(\rho_n\rho_{n+1})^2$
and so
$$
(\rho_n^2)'=2\rho_n\rho_{n+1}(\rho_n\rho_{n+1})'+(\ov\a_{n+1}\rho_n)'(\a_{n+1}\rho_n)+
(\ov\a_{n+1}\rho_n)(\a_{n+1}\rho_n)'. $$
 Hence
\begin{equation}\label{2.30}
 \rho_n'=\rho_{n+1}(\rho_n\rho_{n+1})'+\Re[\a_{n+1}(\ov\a_{n+1}\rho_n)'].
\end{equation}
Next, $(\a_{n-1}\rho_n)'=\a_{n-1}'\rho_n+\a_{n-1}\rho_n'$, so that
\begin{equation}\label{2.31}
 \a_{n-1}'=\frac1{\rho_n}\[(\a_{n-1}\rho_n)'-\a_{n-1}\rho_n'\].
\end{equation}

 The right hand side of (\ref{2.31}) can be expressed in terms of
derivatives of the CMV matrix entries and thereby, via the Lax
equation, of $\a$'s themselves. First, by (\ref{2.20}) and
(\ref{2.23})
\begin{eqnarray}\label{2.32}
  (\rho_n\rho_{n+1})' &=& c_{n,n+2}'+c_{n+2,n}'=K_{n,n+2}+K_{n+2,n} \\
  {} &=& \rho_n\rho_{n+1}\Re(\ov\a_n\a_{n-1}-\ov\a_{n+2}\a_{n+1}).
\end{eqnarray}
Next, it is immediate from (2.3), (2.4) that
$\ov\a_{n+1}\rho_n=c_{n+1,n}\e_n+c_{n,n+1}\e_{n+1}$, and so by
(\ref{2.21}), (\ref{2.22})
\begin{eqnarray}\label{2.33}
  (\ov\a_{n+1}\rho_n)' &=& c_{n+1,n}'\e_n+c_{n,n+1}'\e_{n+1}=
  K_{n+1,n}\e_n+K_{n,n+1}\e_{n+1}\\
  {} &=&
  \ov\a_{n+1}\rho_n\Re(\ov\a_n\a_{n-1}-\ov\a_{n+1}\a_n)+\rho_n\rho_{n+1}^2\ov\D_{n+1}.
\end{eqnarray}
Similarly
$$ (\a_{n-1}\rho_n)' =
   \a_{n-1}\rho_n\Re(\ov\a_n\a_{n-1}-\ov\a_{n+1}\a_n)+\rho_n\rho_{n-1}^2\D_{n-1}.$$
Plugging (\ref{2.32}) and (\ref{2.33}) into (\ref{2.30}) gives
\begin{eqnarray*}
  \rho_n' &=& \rho_n\rho_{n+1}^2\Re(\ov\a_n\a_{n-1}-\ov\a_{n+2}\a_{n+1}
  +\a_{n+1}\ov\D_{n+1})
  +\rho_n|\a_{n+1}|^2\Re(\ov\a_n\a_{n-1}-\ov\a_{n+1}\a_{n}) \\
  {} &=& \rho_n\rho_{n+1}^2\Re(-\ov\a_{n+2}\a_{n+1}+\ov\a_{n+1}\a_{n}
  +\a_{n+1}\ov\D_{n+1})+
  \rho_n\Re(\ov\a_n\a_{n-1}-\ov\a_{n+1}\a_{n}) \\
  {} &=& \rho_n\Re(\ov\a_n\a_{n-1}-\ov\a_{n+1}\a_{n}).
  \end{eqnarray*}
In the upshot, the Schur flow equations emerge from (\ref{2.31}):
 \begin{eqnarray*}
\a_{n-1}'\rho_n &=&
\a_{n-1}\rho_n\Re(\ov\a_n\a_{n-1}-\ov\a_{n+1}\a_n)+\rho_n\rho_{n-1}^2\D_{n-1} \\
{} &=&
-\a_{n-1}\rho_n\Re(\ov\a_n\a_{n-1}-\ov\a_{n+1}\a_n)=\rho_n\rho_{n-1}^2\D_{n-1},
\end{eqnarray*}
as claimed. \hfill $\square$

\noindent
{\sl Remark}. We could equally well have considered the
general IBV problem, that, strictly speaking has nothing to do
with OPUC:
\begin{equation}\label{}
    \nonumber \a_n'(t)=(1-|\a_n(t)|^2)(\a_{n+1}(t)-\a_{n-1}(t)),
    \qquad t>0
\end{equation}
with a continuous boundary function $|\a_{-1}(t)|\leq 1$. The
above evaluation shows that the Lax form of such problem is
$\mathcal{C}_g'=[A_g,\mathcal{C}_g]$ with
\begin{equation}\label{}
\nonumber \mathcal{C}_g=\left(%
\begin{array}{cccccc}
  -\ov\a_0\a_{-1} & \ov\a_1\rho_0 & \rho_0\rho_1 & 0 & 0 & \ldots \\
  -\a_{-1}\rho_0 & -\ov\a_1\a_0 & -\a_0\rho_1 & 0 & 0 & \ldots \\
  0 & \ov\a_2\rho_1 & -\ov\a_2\a_1 & \ov\a_3\rho_2 & \rho_2\rho_3 & \ldots \\
  0 &\rho_1\rho_2 & -\a_1\rho_2 & -\ov\a_3\a_2 & -\a_2\rho_3 & \ldots \\
  0 & 0 & 0 & \ov\a_4\rho_3 & -\ov\a_4\a_3 & \ldots \\
  \vdots & \vdots & \vdots & \vdots & \vdots & \vdots \\
\end{array}%
\right) \end{equation}
 and
\begin{equation}\label{}
\nonumber A_g=\left(%
\begin{array}{cccccc}
  -\Re\ov\a_0\a_{-1} & \rho_0\ov\D_0 & \rho_0\rho_1 & {} & {} & {} \\
  {} & -\Re\ov\a_1\a_0 & \rho_1\D_1 & \rho_1\rho_2 & {} & {} \\
  {} & {} & -\Re\ov\a_2\a_1 & \rho_2\ov\D_2 & \rho_2\rho_3 & {} \\
  {} & {} & {} & \ddots & \ddots & \ddots \\
\end{array}%
\right)  \end{equation} (to be compared with (\ref{1.17}) and
(\ref{1.22})). Furthermore, the doubly infinite system
\begin{equation}\label{}
    \nonumber \a_n'(t)=(1-|\a_n(t)|^2)(\a_{n+1}(t)-\a_{n-1}(t)),
    \qquad t>0, \quad n\in\mathbb{Z}
\end{equation}
with the initial data
$\{\a_n(0)\}_{k\in\mathbb{Z}}\subset\mathbb{D}$ has its equivalent
Lax form
$$ \hat{\mathcal{C}}'=[\hat A,\hat{\mathcal{C}}], $$
where $\hat{\mathcal{C}}$ and $\hat A$ are doubly infinite
extensions of (\ref{1.17}) and (\ref{1.22}) given by the same
expressions (2.1)--(2.4) and (2.5)--(2.6), respectively, with
$$ \e_m=\frac{1-(-1)^{|m|}}2\,, \qquad m\in\mathbb{Z}. $$

\section{Proof of Theorem 1: $(2)\Rightarrow (3)\Rightarrow (1)$.}

\noindent
$(2)\Rightarrow (3)$. Let
 $$ R_z(t):= (\mathcal{C}-zI)^{-1}=\|r_{mn}(z,t)\|_{m,n=0}^\infty
 $$
 be a resolvent of the CMV matrix (\ref{1.17}). It is easy to see
 that $R_z$ obeys the same Lax equation (\ref{1.21}). Indeed,
 differentiating the identity $(\mathcal{C}-zI)R_z=I$ with respect
 to $t$ entails
 $$ \mathcal{C}'(t)R_z(t)+(\mathcal{C}(t)-zI)R_z'(t)=0, $$
 and so
 \begin{equation}\label{3.1}
  R_z'(t)=-R_z(t)\mathcal{C}'(t)R_z(t)=-R_z(t)[A,\mathcal{C}]R_z(t)=[A,R_z],
\end{equation}
as claimed.

Take the equation for $(0,0)$-entry of (\ref{3.1}):
 \begin{eqnarray*}
 r_{00}' &=& \Re\ov\a_0r_{00}+\rho_0\ov\D_0 r_{10}+\rho_0\rho_1
 r_{20}-r_{00}\Re\ov\a_0 \\
 &=& \rho_0(\ov\a_1+1)r_{10}+\rho_0\rho_1 r_{20}.
 \end{eqnarray*}
As it follows from $(\mathcal{C}-zI)R_z=I$
 $$ (\ov\a_0-z)r_{00}+\ov\a_1\rho_0 r_{10}+\rho_0\rho_1 r_{20}=1,
 $$
which allows to eliminate $r_{20}$ in favor of $r_{10}$, $r_{00}$,
and so
\begin{equation}\label{3.5}
 r_{00}'=\rho_0 r_{10}+1-(\ov\a_0-z)r_{00}.
\end{equation}

By the Spectral Theorem, the resolvent entries can be found from
$$
r_{m,n}(z,t)=\int_\mathbb{T}\frac{\chi_n(\z)\ov{\chi_m(\z)}}{\z-z}\,
d\mu(\z,t) $$
 with $\chi_n$ (\ref{1.14}), which is particularly  simple for the first
 two elements
$$ r_{00}(z,t)=\int_{\mathbb{T}} \frac{d\mu(\z,t)}{\z-z}\,, \qquad
   r_{10}(z,t)=\int_{\mathbb{T}}
   \frac{\ov{\f_1(\z)}}{\z-z}\,d\mu(\z,t),
$$
 $\f_1(z)=\rho_0^{-1}(z-\ov\a_0)$. Hence for the right hand
 side of (\ref{3.5})
$$ \rho_0 r_{10}+1-(\ov\a_0-z)r_{00}=
   \int_{\mathbb{T}}\frac{\z+\ov\z-2\Re\a_0}{\z-z}\,d\mu(\z,t) $$
holds, and we end up with a differential equation for
orthogonality measures
$$ d\mu'(\z,t)=(\z+\ov\z-2\Re\a_0)d\mu(\z,t). $$
Finally,
\begin{eqnarray*}
d\mu(\z,t) &=&
\exp\left\{\int_0^t(\z+\ov\z-2\Re\a_0(s))ds\right\}\,d\mu(\z,0)\\
{} &=&  C(t)e^{t(\z+\z^{-1})}\,d\mu(\z,0),
\end{eqnarray*}
as needed.

\noindent
$(3)\Rightarrow (1)$. We start out from the
transformation of the orthogonality measure
$d\mu(\z,t)=\f(\z,t)\,d\mu(\z,0)$ and derive a differential
equation for the moments
$$ \mu_k(t)=\int_{\mathbb{T}} \z^{-k}\,d\mu(\z,t)=
   \int_{\mathbb{T}} \z^{-k}\f(\z,t)\,d\mu(\z,0). $$
Now
$$ \f(\z,t)=\frac{\exp(t(\z+\z^{-1}))}{f(t)}, \qquad f(t)=
   \int_{\mathbb{T}} e^{t(\z+\z^{-1})}\,d\mu(\z,0), $$
   and so
$$
\f'(\z,t)=(\z+\z^{-1})\f(\z,t)-\frac{f'(t)}{f(t)}\f(\z,t).
$$
As $f'f^{-1}=c_1+c_{-1}$, we have
\begin{equation}\label{3.10}
 \mu'_k(t)=\mu_{k-1}(t)+\mu_{k+1}(t) - g(t)\mu_k(t),
 \qquad g(t)=c_1(t)+c_{-1}(t),
 \qquad k\in\mathbb{Z}.
\end{equation}

The rest is based heavily on (\ref{1.12}) which relates Verblunsky
coefficients and moments of the orthogonality measure. The idea to
differentiate determinants and take into account (\ref{3.10}) goes
back to \cite{Be1}, see also \cite [Lemma 1]{MN1}. For a set of
integers $k_1<k_2<\ldots<k_n$ denote
$$ T(k_1,\ldots,k_n):=
\left(%
\begin{array}{cccc}
  \mu_{k_1} & \mu_{k_1-1} & \ldots & \mu_{k_1-n+1} \\
  \mu_{k_2} & \mu_{k_2-1} & \ldots & \mu_{k_2-n+1} \\
  \vdots & \vdots & \vdots & \vdots \\
  \mu_{k_n} & \mu_{k_n-1} & \ldots & \mu_{k_n-n+1} \\
\end{array}%
\right), $$ $D(k_1,\ldots,k_n):= \det T(k_1,\ldots,k_n)$, and so
$D_n=D(0,1,\ldots,n-1)$. Put $G_n=D(-1,0,\ldots,n-2)$ and write
(\ref{1.12}) as $\F_n(0,t)=(-1)^nG_n D_n^{-1}$. Then
 \begin{equation}\label{3.12}
  \F_n'(0,t)=(-1)^n\frac{G_n'D_n-G_nD_n'}{D_n^2}. \end{equation}
  It is clear from (\ref{3.10}) that intermediate determinants in the sum
$$ D_n'=
\left|%
\begin{array}{cccc}
  \mu_0' & \mu_{-1}' & \ldots & \mu_{-n+1}' \\
  \mu_1 & \mu_0 & \ldots & \mu_{-n+2} \\
  \vdots & \vdots & \vdots & \vdots \\
  \mu_{n-1} & \mu_{n-2} & \ldots & \mu_0 \\
\end{array}%
\right| + \ldots +
\left|%
\begin{array}{cccc}
  \mu_0 & \mu_{-1} & \ldots & \mu_{-n+1} \\
  \mu_1 & \mu_0 & \ldots & \mu_{-n+2} \\
  \vdots & \vdots & \vdots & \vdots \\
  \mu_{n-1}' & \mu_{n-2}' & \ldots & \mu_0' \\
\end{array}%
\right| $$ have the same value $-gD_n$, whereas the first and the
last ones equal, respectively,
\begin{eqnarray*}
\left|%
\begin{array}{cccc}
  \mu_0' & \mu_{-1}' & \ldots & \mu_{-n+1}' \\
  \mu_1 & \mu_0 & \ldots & \mu_{-n+2} \\
  \vdots & \vdots & \vdots & \vdots \\
  \mu_{n-1} & \mu_{n-2} & \ldots & \mu_0 \\
\end{array}%
\right| &=& D(-1,1,2,\ldots,n-1)-g(t)D_n, \\
\medskip
\left|%
\begin{array}{cccc}
  \mu_0 & \mu_{-1} & \ldots & \mu_{-n+1} \\
  \mu_1 & \mu_0 & \ldots & \mu_{-n+2} \\
  \vdots & \vdots & \vdots & \vdots \\
  \mu_{n-1}' & \mu_{n-2}' & \ldots & \mu_0' \\
\end{array}%
\right| &=& D(0,1,\ldots,n-2,n)-g(t)D_n.
\end{eqnarray*}
Hence
\begin{equation}\label{3.14}
D_n'(t)=-ng(t)D_n(t)+D(-1,1,\ldots,n-2,n-1)+D(0,1,\ldots,n-2,n).
\end{equation}
Similarly,
\begin{equation}\label{3.15}
G_n'(t)=-ng(t)G_n(t)+D(-2,0,\ldots,n-3,n-2)+D(-1,0,\ldots,n-3,n-1).
\end{equation}
After plugging (\ref{3.14}) and (\ref{3.15}) into (\ref{3.12}) we
come to
\begin{eqnarray}\label{3.20}
\nonumber  \F_n'(0,t) &=& \frac{(-1)^n}{D_n^2}\left\{
D_n\[D(-2,0,\ldots,n-2)+D(-1,\ldots,n-3,n-1)\]\right. \\
\nonumber &-& \left.
G_n\[D(-1,1,\ldots,n-1)+D(0,\ldots,n-2,n)\]\right\}.
\end{eqnarray}

Let us now go over to the right hand side of (\ref{1.20}), written
for $\F_n(0)$:
$$ \(1-|\F_n(0)|^2\)\(\F_{n+1}(0)-\F_{n-1}(0)\)=
   (-1)^{n+1}\frac{D_n^2-|G_n|^2}{D_n^2}
   \left\{\frac{G_{n+1}}{D_{n+1}}-\frac{G_{n-1}}{D_{n-1}}\right\}.
$$
The standard Silvester identity applied to the matrix
$T(0,1,\ldots,n)$ gives
$$ D_n^2-|G_n|^2=D_{n+1}D_{n-1}, $$
and so
\begin{equation}\label{}
\nonumber \(1-|\F_n(0)|^2\)\(\F_{n+1}(0)-\F_{n-1}(0)\)=
   \frac{(-1)^{n+1}}{D_n^2}(G_{n+1}D_{n-1}-D_{n+1}G_{n-1}).
\end{equation}
Another application of the Silvester identity (in a bit modified
form) shows that
\begin{eqnarray*}
  G_{n+1}D_{n-1} &=& G_n D(-1,1,\ldots,n-1)-D_n D(-2,0,\ldots,n-2) \\
  D_{n+1}G_{n-1} &=& D_n D(-1,0,\ldots,n-3,n-1)-G_n
  D(0,1,\ldots,n-2,n),
\end{eqnarray*}
and we arrive at the Schur flow (\ref{1.20}). That completes the
proof of Theorem 1. \hfill $\square$

\smallskip

There is yet another way to prove $(3)\Rightarrow (1)$, which
gives not only (\ref{1.20}), but the differential equations for
the monic orthogonal polynomials. I learned it from \cite [Section
8.3]{Is}.

\noindent
{\bf Theorem 2}. {\it The monic polynomials
$\F_n(\cdot,t)$ orthogonal with respect to the Bessel
transformation $(\ref{1.23})$ satisfy the differential equation}
\begin{equation}\label{3.50}
 \F_n'(z,t)=\F_{n+1}(z,t)-(z+\ov\a_{n}\a_{n-1})\F_n(z,t)-
 (1-|\a_{n-1}|^2)\F_{n-1}(z,t).
\end{equation}

\noindent
{\it Proof}. The idea is to differentiate orthogonality
relations (\ref{1.7}) with respect to $t$. Now $w'=(\z+\ov\z)w$,
and we have for $m<n$
 \begin{equation}\label{3.52}
   \int_{\mathbb{T}} \F_m'\ov\F_n\,d\mu(\z,t) +
   \int_{\mathbb{T}} \F_m\ov{\F_n'}\,d\mu(\z,t) +
   \int_{\mathbb{T}} (\z+\ov\z)\F_m'\ov\F_n\,d\mu(\z,t)=0.
   \end{equation}
It is clear that $\F_m'$ is a polynomial of degree at most $m-1$,
and so the first integral in the above sum is zero. If $m\leq n-2$
then
$$ \int_{\mathbb{T}} \F_m\ov{(\F_n'+\z\F_n)}\,d\mu(\z,t)=0, $$
and hence
 \begin{equation}\label{3.54}
 \F_n'(z)+z\F_n(z)=x_n\F_{n+1}(z)+y_n\F_n(z)+z_n\F_{n-1}(z)
  \end{equation}
  with some parameters $x_n$, $y_n$, $z_n$ depending on $t$. By
matching the coefficients for $z^{n+1}$ and $z^n$ in (\ref{3.54})
and using (\ref{1.6}) we find
 $$ x_n=1, \qquad y_n=l_n(t)-l_{n+1}(t). $$
To get $z_n$ take (\ref{3.52}) with $m=n-1$ and apply the \SZ\
recurrences (\ref{1.10}):
$$ 0=\int_{\mathbb{T}} \F_{n-1}\ov{(\F_n'+\z\F_n)}\,d\mu(\z,t)+
     \int_{\mathbb{T}}
     (\F_n+\ov\a_{n-1}\F_{n-1}^*)\ov\F_n\,d\mu(\z,t), $$
and so
\begin{eqnarray*}
  \int_{\mathbb{T}} \F_{n-1}\ov{(\F_n'+\z\F_n)}\,d\mu(\z,t) &=& -\,\|\F_n\|^2_{\mu}, \\
  z_n &=&
  -\,\frac{\|\F_n\|^2_{\mu}}{\|\F_{n-1}\|^2_{\mu}}=
  -\frac{\k_{n-1}^2}{\k_n^2}=|\a_{n-1}|^2-1.
\end{eqnarray*}

To find the expression for $l_n$ we turn to (\ref{1.18})
$$
\F_n(z)=z^n+l_nz^{n-1}+\ldots=\prod_{j=0}^{n-1}(z+\ov\a_j\a_{j-1})+\ldots,
$$
so that $l_n=\sum_{j=0}^{n-1}\ov\a_j\a_{j-1}$, and we come to
(\ref{3.50}). Putting $z=0$ yields (\ref{1.20}), as was to be
proved. \hfill $\square$

\smallskip

It might be worth pointing out that some properties of Verblunsky
coefficients for the Bessel transformed measures (such as the rate
of decay) are inherited from those of the initial data.

\noindent
{\bf Theorem 3}.  {\it Let $\a_n(0)$ enjoy either of the
properties
\begin{enumerate}
    \item  $\{\a_n(0)\}\in\ell^p$,\ \  $p=1,2$;
    \item  $|\a_n(0)|\leq C e^{-\a n}$,\ \  $\a>0$.
\end{enumerate}
Then the same holds for $\a_n(t)$ for each $t>0$}.

\noindent
{\it Proof}. It is obvious from (\ref{1.23}) that
$\mu(\z,t)$ belongs to the \SZ\ class (i.e.,
$\log\mu'\in\L^1(\mathbb{T}$) if and only if $\mu(\z,0)$ does, and
so the first statement with $p=2$ follows from fundamental Szeg\H
o's Theorem \cite [Theorem 2.3.1]{S1}. As for the case $p=1$, note
that by Baxter's theorem (see, e.g., \cite [Theorem 5.2.1]{S1})
$\mu(\z,0)=wdm$ with $w>0$ and $w\in W$, class of absolutely
convergent Fourier series. It is clear from (\ref{1.23}) and the
Wiener--Levy theorem that
$$ d\mu(\z,t)=w(\z,t)dm, \qquad w(\z,t)=e^{t(\z+\z^{-1})}w(\z)\in
W, $$
 $w(\z,t)>0$, and so the repeated application of Baxter's theorem
 does the job.

 To prove the second statement, we introduce the \SZ\ function
 $$ D(z,\mu):=\exp\(\frac12\int_\mathbb{T} \frac{\z+z}{\z-z}\,\log
 w(\z)\,dm\), $$
 defined for an arbitrary measure $\mu=wdm+\mu_s$ from the \SZ\
 class. A straightforward computation gives for the Bessel transformed
 measures (\ref{1.23})
\begin{equation}\label{3.57}
D(z,\mu(t))=\exp\(\frac12\int_\mathbb{T}\frac{\z+z}{\z-z}\,
  [t(\z+\z^{-1})+\log w(\z)]\,dm\)=e^{tz}D(z,\mu(0)).
  \end{equation}
Denote by $R_t$ the radius of convergence of the Taylor series for
$D^{-1}(z,\mu(t))$ about the origin. By the Nevai--Totik theorem
\cite{NT}
$$ \limsup_{n\to\infty}|\a_n(t)|^{\frac1n}=R_t^{-1}, $$
and it is clear from (\ref{3.57}) that $R_t=R_0$ for all $t>0$, as
claimed. \hfill $\square$

Note that under assumptions of Theorem 3 the series
$$ \mathcal{K}=-\sum_{j=0}^\infty \ov\a_j(t)\a_{j-1}(t) $$
converges absolutely and the Schur flow can be written in the form
$$ \a_j'=\{\a_j,H\}, \qquad H=-2\Im\mathcal{K}, $$
where the (formal) Poisson brackets are defined by
$$ \{f,g\}=i\sum_{j\geq 0} \rho_j^2\[\frac{\partial f}{\partial\ov\a_j}
\frac{\partial g}{\partial\a_j}-\frac{\partial
f}{\partial\a_j}\frac{\partial g}{\partial\ov\a_j}\]. $$
 So (\ref{1.20}) is the evolution of the Verblunsky coefficients
 under the flow generated by the Hamiltonian $-2\Im\mathcal{K}$
 (cf. \cite{Ne1}).

\section{Modified Bessel measures on the unit circle}

Because of the boundary condition $\a_{-1}=-1$ IBV problem
(\ref{1.33}) with zero initial conditions
$$ \a_0(0)=\a_1(0)=\ldots=0 $$
has a nontrivial solution. Theorem 1, (3), shows that we are
dealing now with the Bessel transformation of the Lebesgue measure
$$ d\mu(\z,t)=C(t)e^{t(\z+\z^{-1})}\,dm, $$
called in the sequel the {\it modified Bessel measures on the unit
circle}, with $\b_n(t)$ the Verblunsky coefficients of
$\mu(\cdot,t)$. The corresponding system of orthogonal polynomials
has arisen from studies of the length of longest increasing
subsequences of random words \cite{BDJ} and matrix models
\cite{PSh} (see \cite [example 8.3.4]{Is} for more detail about
the modified Bessel OPUC).

Note first that $C(t)$ can be easily computed
\begin{eqnarray*}
  C^{-1}(t) &=& \int_\mathbb{T} e^{t(\z+\z^{-1})}=
  \frac1{2\pi}\int_0^{2\pi} e^{2t\cos x}\,dx=
  \frac1{2\pi}\sum_{n=0}^\infty \frac{(2t)^n}{n!}\int_0^{2\pi}(\cos x)^n\,dx  \\
  {} &=& \sum_{n=0}^\infty \frac{t^{2n}}{(n!)^2}=I_0(2t),
\end{eqnarray*}
where $I_k$ is the modified Bessel function of order $k$.
Similarly, for the moments of the measure we have
$$ \mu_p(t)=\int_\mathbb{T} \z^{-p}\,d\mu(\z,t)=\frac{I_p(2t)}{I_0(2t)}, \quad
p\in\mathbb{Z}_+, \qquad \mu_{-p}=\mu_p. $$

There is an important feature of the modified Bessel measures,
namely, their Verblunsky coefficients satisfy a nonlinear
recurrence relation known as the discrete Painlev\'e II equation
(see \cite [lemma 8.3.5]{Is})
 \footnote{In the notation of \cite{Is} $r_n(t)=-\b_{n-1}(t/2)$.}
\begin{equation}\label{4.1}
 -(n+1)\frac{\b_n(t)}{t(1-\b_n^2(t))}=\b_{n+1}(t)+\b_{n-1}(t),
 \quad n\in\mathbb{Z}_+,
\end{equation}
with $\b_{-1}=-1$, $\b_0=I_1(2t)/I_0(2t)$. Clearly, all $\b$'s are
real now.

Before studying the long time behavior of $\b_n$ we prove an
auxiliary result concerning the general Schur flows.

\noindent
{\bf Lemma 4.} {\it For the solution $\a_n$ of the Schur
flow $(\ref{1.20})$ the limit relations hold
\begin{equation}\label{4.5}
\lim_{t\to\infty} \rho_n(t)\rho_{n+1}(t)=0, \qquad
\lim_{t\to\infty} \rho_n(t)\D_n(t)=0, \quad n\in\mathbb{Z}_+
\end{equation}
and
\begin{equation}\label{4.7}
\lim_{t\to\infty} \Re\ov{\a_{n}(t)}\a_{n-1}(t)=\g_n, \qquad
n\in\mathbb{Z}_+.
\end{equation}
Furthermore, $\g_n$ is monotonically increasing:
$-1\le\g_0\le\g_1\le\ldots\le 1$. }

\noindent
{\it Proof}. Let us focus on the form of the Lax
 equation given in (\ref{1.35})--(\ref{1.37}). As $B^*=-B$, it is
 clear that
 \begin{equation}\label{4.10}
   L'=[\pi(L),L], \qquad L:=\frac{\mathcal{C}+\mathcal{C}^*}2\,,
   \quad \pi(L)=L_+-L_-.
 \end{equation}
The latter is exactly what is called in \cite{DLT} the {\it Toda
flow}. By \cite [Proposition 5]{DLT} $L(t)$ converges strongly to
a diagonal operator $\mbox{diag}\ (d_0,d_1,\ldots)$, which implies
(\ref{4.5})--(\ref{4.7}) in view of the explicit expression for
$L$. Note that (\ref{4.7}) with $\g_n=-d_n$ comes from the
diagonal entries of $L$, whereas both relations in (\ref{4.5})
from the off diagonal entries.

The second statement follows from (\ref{1.31})
$$
\rho_n(t)=\rho_n(0)\exp\left\{\int_0^t
\Re(\ov\a_{n}(s)\a_{n-1}(s)-\ov\a_{n+1}(s)\a_{n}(s))ds\right\},
$$
and since the left-hand side is bounded as $t\to\infty$ we have
$\g_{n}\le\g_{n+1}$. \hfill $\square$

Let us go back to the modified Bessel measures and their
Verblunsky coefficients, and prove

\noindent
{\bf Theorem 5}. {\it The limit relations
\begin{equation}\label{4.20}
  \lim_{t\to\infty} \b_n(t)=(-1)^n,
\end{equation}
\begin{equation}\label{4.21}
  \lim_{t\to\infty} t(1-\b_n^2(t))=\frac{n+1}2
\end{equation}
hold for all $n\in\mathbb{Z}_+$.}

\noindent
{\it Proof}. Note first that (\ref{4.21}) is an
immediate consequence of (\ref{4.20}) and (\ref{4.1}). To prove
(\ref{4.20}) we proceed in two steps.

\noindent
 1. As we know (see Lemma 4) the sequence
 $\g_n\uparrow$. Let us show that in fact all $\g_n$'s are the same.
 Assume on the contrary that $\g_k<\g_{k+1}$ for some $k$. It
 follows from (\ref{1.31}) that
 $$ \rho_k^2(t)\leq Ce^{-\d t}, \qquad \d>0, \qquad \rho_k^2=1-\b_k^2. $$
 But then
 $$ \frac{|\b_n(t)|}{t(1-\b_k^2(t))}=
    \frac{\sqrt{1-\rho_k^2(t)}}{t\rho_k^2(t)}\ge
    \frac{\sqrt{1-\rho_k^2(t)}}{Ct}\,e^{\d t} \to+\infty, \quad
    t\to\infty $$
that contradicts (\ref{4.1}), since the right-hand side there is
bounded. So we need only find the common value of $\g_n$.

\noindent
 2. We show that $\g_n=-1$ by computing $\lim_{t\to\infty}\b_0(t)$.
 As is well known,
 $$ I_k(t)=\frac{e^t}{\sqrt{2\pi t}}\(1+O\(\frac1t\)\), \quad
    t\to\infty, $$
and so
$$
\lim_{t\to\infty}\b_0(t)=\lim_{t\to\infty}\frac{I_1(2t)}{I_0(2t)}=1,
$$
that is, $\g_0=-1$, as claimed. The desired result follows from
(\ref{4.7}) by induction. \hfill $\square$

\noindent
{\sl Remark}. It might be a challenging problem to give
a direct proof of (\ref{4.20}) based on the explicit formula
(\ref{1.12}), which now takes on the form
$$ \b_n(t)=(-1)^n\frac{\det\|I_{k-j-1}(2t)\|_{0\leq k,j\leq n}}
   {\det\|I_{k-j}(2t)\|_{0\leq k,j\leq n}}\,,\qquad n\in\mathbb{Z}_+, $$
and the complete asymptotic series expansion for the modified
Bessel function (see, e.g., \cite [Chapter 7.8]{Ol})
$$ I_k(t)\simeq \frac{e^t}{\sqrt{2\pi t}}\sum_{j=0}^\infty
(-1)^j\frac{(4k^2-1^2)\ldots(4k^2-(2j-1)^2)}{j! (8t)^j}\,, \quad
t\to\infty. $$
 I managed to carry out the computation for $n=1$, and it seems
 like one needs $n+1$ terms of this series for $\b_n$.

\smallskip
{\bf Acknowledgement}. I thank Yu.M. Berezanskii for drawing my
attention to the problem discussed in the paper, and M. Ismail for
giving a chance to get acquainted with the manuscript of his
ongoing book \cite{Is}. The work was partially supported by INTAS
Research Network NeCCA 03-51-6637 and NATO Collaborative linkage
grant PST. CLG. 979738.

\end{document}